\documentclass[12pt]{article}
\usepackage{amsfonts,amsmath,amsxtra}

\newtheorem{theorem}{Theorem}
\newtheorem{corollary}{Corollary}
\newtheorem{lemma}{Lemma}

\newenvironment{definition}
{\smallskip\noindent{\bf Definition\/}:}{\smallskip\par}
\newenvironment{proposition}
{\smallskip\noindent{\bf Proposition\/}.}{\smallskip\par}
\newenvironment{propositions}
{\smallskip\noindent{\bf Propositions\/}.}{\smallskip\par}

\newenvironment{remark}
{\smallskip\noindent{\bf Remark\/}.}{\smallskip\par}

\newenvironment{proof}
{\noindent{\bf Proof\/}.}{{ $\Box$}\smallskip\par}

\title{Symmetries of some motivic integrals.}
\author{E. Gorsky \footnote{Partially supported by the grants RFBR-007-00593, INTAS-05-7805, NSh-4719.2006.1, and the Moebius Contest fellowship for young scientists.}}
\date{}

\begin{document}

\maketitle

\begin{abstract}
We give an explicit formula for the motivic integrals related to the Milnor number over  spaces of parametrised arcs on the plane 
with fixed tangency orders with the axis. These integrals are rational functions of the parameters and the class of the affine line. Using a set of natural recurrence relations between them,  we prove an unexpected invariance property with respect to the simultaneous inversion of the parameters and the class of the affine line. We also discuss a generalization of this system of recurrence relations whose solutions are also symmetric and satisfy additional differential equations. 
\end{abstract}

\section{Introduction}

From the Weyl conjectures it is known that the Poincare duality implies a functional equation for the zeta function of the smooth projective variety. This equation describes the behavior of the zeta function under the inversion of the parameter.
M. Kapranov (\cite{kapr}) and F. Heinloth (\cite{hein})
proved that some generalizations of this functional equation also holds for the Kapranov zeta function
of curves and abelian variety respectively. 

Another symmetry property can be observed in the knot theory. The Alexander polynomial of a curve is invariant under the change $t\to t^{-1}$. The result of a change $q\to q^{-1}$ in the Jones polynomial gives a Jones polynomial of the mirror knot. These properties are also true for their categorifications, Heegard-Floer knot homologies of Oszvath-Szabo (\cite{hfk}) and Kvovanov homologies (\cite{KhoHo}).

A. Campillo, F. Delgado and S. Gusein-Zade gave in (\cite{cdg},\cite{cdg2}) an interpretation of the Alexander polynomial of an algebraic knot as an integral with respect to the the Euler characteristic over the projectivisation over the space of functions. This integral is a specification under the Euler characteristic homomorphism of a corresponding motivic integral over the space of functions (\cite{cdg3}) , which can be also considered as a certain categorification of the Alexander polynomial, but its relation to the Heegard-Floer homologies is completely unknown. 

Nevertheless, a nice test for the motivic integral (over the space of functions or over the space of arcs) to be a natural generalization of the Alexander polynomial is its behavior under the simultaneous inversion of a parameter and the class of the affine line in the Grothendieck ring of varieties.

A conclusion of this (a bit speculative) discussion is that the motivic integrals, which transform in a nice way under the inversion of the parameter and the affine line, may have a nice geometric sense and vice versa.

This paper deals with a  set of motivic integrals related to the Milnor number (\cite{book}) over
spaces of arcs with fixed tangency orders with the axis. We prove that all of them are rational functions of the parameter and the class of the affine line, and write a set of recurrence relations between them. Using this relations, we give an explicit formula for these integrals\newline (Theorem 2), from which we deduce the symmetry property for them. 

Namely, if $\mu$ denotes the Milnor number, $\mathbb{L}$ is the class of the affine line in the Grothendieck ring of varieties and $$G_{k,m}=\int_{Ord x(t)=k,Ord y(t)=m}t^{\mu},$$
then:

a) $G_{k,m}$ are rational functions of $t$ and $\mathbb{L}$, 

b) $G_{k,m}(t^{-1},\mathbb{L}^{-1})=t^{-2(k-1)(m-1)}\mathbb{L}^{2(k+m)}\cdot G_{k,m}(t,\mathbb{L}).$

The proof of these facts uses recurrence relations between $G_{k,m}$ for different $k$ and $m$, and the direct proof  is not known.

 We also discuss in the section 4 a certain deformation of the set of equations. Since they arise from the consideration of the blowup of the plane, the class of the exceptional divisor in the Grothendieck ring
is related with their coefficients. Iterative use of these equations corresponds to the sequence of successive blowups. In the physical literature the volumes (or the "Kaehler classes") of the exceptional divisors sometimes (e.g. \cite{phys}) are considered as independent parameters. Motivated by this, we deform our system of equation in such a way that some of classes of the exceptional divisors become parameters that are not necessary equal to $\mathbb{L}+1$. Despite that the geometric meaning of this procedure is far from understanding, we prove an analogue of the Theorem 2 and the symmetry property analogous to (b). We also prove that the generating function for deformations of the motivic measures of sets of arcs with fixed tangency order with axis, that is, the deformed versions of the above integral at $t=1$, satisfies   
a set of differential-functional equations (Theorem 4).
 
\section{Motivic measure}

Let $\mathcal{L}=\mathcal{L}_{\mathbb{C}^2,0}$ be the space
of arcs at the origin on the plane. 
It is the set of pairs $(x(t),y(t))$ of formal power series 
(without degree 0 terms). 
Let $\mathcal{L}_n$ be the space of $n$-jets of such arcs, let $\pi_{n}:\mathcal{L}\rightarrow \mathcal{L}_{n}$ be the natural projection.

Let $K_0(Var_{\mathbb{C}})$ be the Grothendieck ring of quasiprojective
complex algebraic varieties. It is generated by the isomorphism classes
of complex quasiprojective algebraic varieties modulo relations
 $[X]=[Y]+[X\setminus Y],$ where $Y$ is a Zariski closed subset of  $X$. 
Multiplication is given by the formula $[X]\cdot [Y]=[X\times Y].$
Let $\mathbb{L}\in K_0(Var_{\mathbb{C}})$ be the class of the complex line.

Consider the ring $K_0(Var_{\mathbb{C}})[\mathbb{L}^{-1}]$
with the following filtration: $F_k$ is generated by the elements
of type $[X]\cdot [\mathbb{L}^{-n}]$ with $n-\dim X\ge k$. Let $\mathcal{M}$
be the completion of $K_0(Var_{\mathbb{C}})[\mathbb{L}^{-1}]$ corresponding
to this filtration.

On an algebra of subsets of the space $\mathcal{L}$ J. Denef and F. Loeser (\cite{dl}) 
(after M. Kontsevich) have constructed a measure $\chi_g$ with values in the ring  $\mathcal{M}$. 

A subset $A\subset \mathcal{L}$ is said to be cylindric if there exist  $n$ and a constructible set $A_n\subset \mathcal{L}_{n}$ such that
$A=\pi_{n}^{-1}(A_n)$. For the cylindric set $A$ define
$$\mu(A)=[A_n]\cdot \mathbb{L}^{-2n}.$$
It was proved in \cite{dl}, that this measure can be extended to an additive measure on a suitable algebra of subsets in $\mathcal{L}$.    

A function $f:\mathcal{L}\rightarrow G$ with values in an abelian group $G$ is called simple, if its image is countable or finite, and for every $g\in G$ the set $f^{-1}(g)$ is measurable. Using this measure, one can define in the natural way the (motivic) integral for simple functions on $\mathcal{L}$ as
$\int_{\mathcal{L}}f d\mu=\sum_{g\in G} g\cdot\mu(f^{-1}(g)),$
if  the right hand side sum converges in  $G\otimes \mathcal{M}.$

Let $h:Y\rightarrow X$ be a proper birational morphism
of smooth manifolds of  dimension $d$
and $J=h^{*}K_{X}-K_{Y}$ be the relative canonical divisor
on  $Y$ (locally it is defined by the Jacobi determinant).
It defines a function $ord_{J}$ on the space of arcs on $Y$ --
the intersection number between the arc and the divisor.
Then one has the following change of variables formula in the motivic 
integral:
  
\begin{theorem}(\cite{dl})
Let $A$ be a measurable subset in the space of arcs on
$X$, let $\alpha$ be a simple function. Then
$$\int_{A}\alpha d\chi_g=\int_{h^{-1}(A)}(h^{*}\alpha)\mathbb{L}^{-ord_{J}}d\chi_g.$$
\end{theorem}

If $h$ is a blow-up of the origin in the plane, the relative
canonical divisor coincides with the exceptional line, so the
function $ord_{J}$ coincides with the intersection number with this line.

\section{Calculations}

Let $\mu(\gamma)$ be a Milnor number of a germ of plane curve $\gamma(t)=(x(t),y(t))$.
Consider a motivic integral
$$\int_{\mathcal{L}}t^{\mu(\gamma)}d\gamma.$$
It is useful to consider also integrals of a form
$$G_{k,m}(t)=\int_{\{\mbox{\rm Ord}x(t)=k, \mbox{\rm Ord}y(t)=m\}}t^{\mu(\gamma)}d\gamma.$$

We will use a following fact: if $\widetilde{\gamma}$ is a strict transform of $\gamma$ under the blowup of the plane at the origin, and $v(\gamma)$ is the order of $\gamma$ at the origin, then 
\begin{equation}
\mu(\gamma)=\mu(\widetilde{\gamma})+v(\gamma)(v(\gamma)-1).
\end{equation}

\begin{propositions}
1. $G_{k,m}(t)$ is a rational function of the parameter $t$ and the class of the affine line $\mathbb{L}$.

\begin{equation}
\label{2}
G_{1,1}(t)=(\mathbb{L}-1)^2\cdot \mathbb{L}^{-2}.
\end{equation}
\begin{equation}
\label{3}
G_{k,m}(t)=G_{m,k}(t).
\end{equation}
\begin{equation}
\label{4}
G_{k,m}(t)=t^{k(k-1)}\mathbb{L}^{-k}\cdot G_{k,m-k}(t)\,\,\,\,\,\mbox{\rm for}\,\,\,\,\, m>k
\end{equation}
\begin{equation}
\label{5}
G_{k,k}(t)=(\mathbb{L}-1)\sum_{m>k} G_{k,m}(t).
\end{equation}
\begin{equation}
\label{6}
G_{k,k}(t)=\frac{(\mathbb{L}-1)t^{k(k-1)}\mathbb{L}^{-k}}{1-t^{k(k-1)}\mathbb{L}^{1-k}}\sum_{m<k} G_{k,m}(t).
\end{equation}

\end{propositions}

\begin{proof}
Propositions (2) and (3) are obvious from the definition of $G_{k,m}$. Proposition 4 follows from the formula of transformation of the Milnor number under the blowup mentioned above and the change of variables formula for the motivic integrals (\cite{dl}). If orders of the series
$x(t)$ and $y(t)$ are equal, then $y(t)=\lambda x(t)+y'(t)$, where $\lambda$ is a nonzero number , and the order of  $y'(t)$ is greater than the one of $y(t)$. This implies the proposition (5). Now $$G_{k,k}(t)=(\mathbb{L}-1)\sum_{m>k} G_{k,m}(t)=(\mathbb{L}-1)^2\sum_{m>k}t^{k(k-1)}\mathbb{L}^{-k}G_{k,m-k}=$$

$$=(\mathbb{L}-1)^2\sum_{m>k}t^{k(k-1)}\mathbb{L}^{-k}(\sum_{m<k} G_{k,m}(t)+G_{k,k}(t)+\sum_{m>k} G_{k,m}(t))=$$

$$=(\mathbb{L}-1)^2\sum_{m>k}t^{k(k-1)}\mathbb{L}^{-k}(\sum_{m<k} G_{k,m}(t)+G_{k,k}(t)(1+\frac{1}{\mathbb{L}-1})),$$
so we get a proposition 6. The proposition 1 now follows from 2, 3, 4 and 6.
\end{proof}

\begin{definition}
Let us introduce the following bivariate polynomials:
$$S_{a,k}(t,\mathbb{L}^{-1})=\sum_{1\le m<k, \mbox{\rm g.c.d}(m,k)=a}t^{(k-1)(m-1)-(a-1)^2}\mathbb{L}^{2a-k-m}.$$
\end{definition}

\begin{propositions} 
\begin{equation}
\label{7}
\mbox{\rm If g.c.d$(k,m)=a$, then}\,\,\,\,\,\,\ G_{k,m}(t,\mathbb{L})=t^{(k-1)(m-1)-(a-1)^2}\mathbb{L}^{2a-k-m}G_{a,a}(t,\mathbb{L}).
\end{equation}
\begin{equation}
\label{8} G_{k,k}(t,\mathbb{L})=\frac{(\mathbb{L}-1)t^{k(k-1)}\mathbb{L}^{-k}}{1-t^{k(k-1)}\mathbb{L}^{1-k}}\sum_{a|k} S_{a,k}(t,\mathbb{L}^{-1})\cdot G_{a,a}(t).
\end{equation}
\begin{equation}
\label{9} G_{k,k}(t,\mathbb{L})=\sum_{1=a_0,a_1,\ldots,a_r=k}\frac{(\mathbb{L}-1)^{2+r}t^{\sum_{j=1}^{r}a_{j}(a_j-1)}\mathbb{L}^{-2-\sum_{j=1}^{r}a_j}}{\prod_{j=1}^{r}(1-t^{a_j(a_j-1)}\mathbb{L}^{1-a_j})}\prod_{j=1}^{r}S_{a_{j-1},a_j}(t,\mathbb{L}^{-1}),
\end{equation}
where summation is indexed by all tuples $(1=a_0<a_1<\ldots<a_r=k)$ such that $a_{j-1}|a_j$
for all $j$.
\end{propositions}

\begin{proof}
Note that $(k-1)(m-1)=k(k-1)+k(m-k-1)$, so,
since the Euclid's algorithm for $k$ and $m$ finishes at $a$, property (4) implies (7). 
Substituting (7) in the proposition 6, we get (8).
We apply (7) and (8) until we reach $a=1$, and using (2) we finally get the equation (9).
\end{proof}

Let us calculate now $S_{a,k}(t,\mathbb{L}^{-1}).$ Let (for $a|k$)
$$\hat{S}_{a,k}(t,\mathbb{L}^{-1})=\sum_{1\le m<k, a|m}t^{(k-1)(m-1)}\mathbb{L}^{-k-m}.$$
In the right hand side we have a sum of a geometric series, so
$$\hat{S}_{a,k}(t,\mathbb{L}^{-1})=\frac{t^{(k-1)(a-1)}\mathbb{L}^{-k-a}-t^{(k-1)^2}\mathbb{L}^{-2k}}{1-t^{(k-1)a}\mathbb{L}^{-a}}.$$

From the other hand, $$\hat{S}_{a,k}(t,\mathbb{L}^{-1})=\sum_{a|b, b|k}t^{(b-1)^2}\mathbb{L}^{-2b}S_{b,k}(t,\mathbb{L}^{-1}).$$
The Moebius inversion formula implies the following equation.

\begin{lemma}
$$t^{(a-1)^2}\mathbb{L}^{-2a}S_{a,k}(t,\mathbb{L}^{-1})=\sum_{b\vdots a, b|k}\mu(b/a)\hat{S}_{b,k}(t,\mathbb{L}^{-1}),$$
where $\mu$ -- denote the Moebius function.
\end{lemma}

From this this lemma we immediately get a formula for $S_{a,k}$.

\begin{proposition}
$$S_{a,k}(t,\mathbb{L}^{-1})=t^{-(a-1)^2}\mathbb{L}^{2a}\sum_{b\vdots a, b|k}\mu(\frac{b}{a})\frac{t^{(k-1)(b-1)}\mathbb{L}^{-k-b}-t^{(k-1)^2}\mathbb{L}^{-2k}}{1-t^{(k-1)b}\mathbb{L}^{-b}}.$$
\end{proposition}

Let us substitute it to the equation (9) and slightly simplify it. We get a following statement.

\begin{theorem}
Let $a=$g.c.d.$(k,m)$. Then
$$G_{k,m}(t,\mathbb{L})=\int_{\{\mbox{\rm Ord}x(t)=k, \mbox{\rm Ord}y(t)=m\}}t^{\mu(\gamma)}d\gamma=$$
$$=(\mathbb{L}-1)^2\cdot t^{(k-1)(m-1)}\mathbb{L}^{-k-m}\cdot \sum_{\underline{a},\underline{b}}\prod_{j=1}^{r}\mu(\frac{b_j}{a_{j-1}})\cdot \frac{(\mathbb{L}-1)t^{(a_j-1)b_j}\mathbb{L}^{-b_j}(1-t^{(a_j-1)(a_j-b_j)}\mathbb{L}^{b_j-a_j})}{(1-t^{a_j(a_j-1)}\mathbb{L}^{1-a_j})\cdot (1-t^{(a_j-1)b_j}\mathbb{L}^{-b_j})},$$
where summation is indexed by all tuples $$(1=a_0\le b_1<a_1\ldots a_{r-1}\le b_r<a_r=a)$$ such that  $a_{j-1}|b_j$ and $b_j|a_j$
for all $j$ ( $r$ is not fixed).

Recall that in the left hand side $\mu$ denotes the Milnor number, and in the right hand side -- the Moebius function, but these appearances of the same symbol in two different meanings should not confuse.
\end{theorem}

\begin{corollary}
$$G_{k,m}(t^{-1},\mathbb{L}^{-1})=t^{-2(k-1)(m-1)}\mathbb{L}^{2k+2m-2}\cdot G_{k,m}(t,\mathbb{L}).$$
\end{corollary}

\begin{proof}
This equation follows from the invariance of every of  $r$ multipliers in the right hand side of the theorem 1 under the change $(t,\mathbb{L})\mapsto (t^{-1},\mathbb{L}^{-1}).$  
\end{proof}

Consider a motivic integral 
$$F(t,\mathbb{L};a,b,c,d,e)=\int_{\mathcal{L}}t^{\mu}a^{v_{x}}b^{v_{y}}c^{v_{x}^2}d^{v_{x}v_y}e^{v_{y}^2}d\mu=\sum_{k,m=1}^{\infty}G_{k,m}(t,\mathbb{L})a^{k}b^{m}c^{k^2}d^{km}e^{m^2}.$$
Analogously to the calculations made above we can write a functional equation
$$F(t,a,b,c,d,e)=F(t, t^{-1}ab\mathbb{L}^{-1}, b, tcde, de^2, e)+
F(t, t^{-1}ab\mathbb{L}^{-1}, a, tcde, dc^2, c)+$$
$$+F(t, t^{-1}ab\mathbb{L}^{-1}, 1, tcde, 1, 1)\cdot(\mathbb{L}-1).$$
It is shown in \cite{eqneng}, that its solution can be completely reconstructed from the initial condition (2) .

The corollary 1 shows that every solution of this functional equation  (in the class of formal power series of $a,b,c,d,e$) satisfies an unexpected symmetry property

$$F(t,\mathbb{L};a,b,c,d,e)=t^2\mathbb{L}^{2}F(t^{-1},\mathbb{L}^{-1};at^{-2}\mathbb{L}^{-2},bt^{-2}\mathbb{L}^{-2},c,dt^{2},e).$$

\section{Additional parameters}

We propose a following generalization of the set of recurrence relations on $G_{k,m}$:
\begin{equation}
\label{lambda}
\begin{cases}
&G_{k,m}(t)=G_{m,k}(t)\cr

&G_{k,m}(t)=t^{k(k-1)}\mathbb{L}^{-k}\cdot G_{k,m-k}(t)\,\,\,\, \mbox{\rm for}\,\,\,\,\, m>k \cr

&G_{k,k}(t)=(\lambda_k-1)\sum_{m>k} G_{k,m}(t)\,\,\,\,\,\,\,

\end{cases}
\end{equation}
where $\lambda_k$ are additional parameters.

\begin{remark} We have $G_{1,k}=\mathbb{L}^{1-k}G_{1,1},$ so $\sum_{k>1} G_{1,k}=\frac{\mathbb{L}^{-1}}{ 1-\mathbb{L}^{-1}}G_{1,1}=\frac{1}{ \mathbb{L}-1}G_{1,1}.$ Therefore $\lambda_1=\mathbb{L},$
if (\ref{lambda}) is not controversial. 
\end{remark}

This choice of parameters is not completely occasional and has some motivation from the mathematical physics. At first, we can remark that  appearance of $\mathbb{L}$ in the equations (4) and (5)  has different nature: equation (4) arise from the change of variables formula, and $\mathbb{L}$ there comes from this formula, so it is related to the role of $\mathbb{L}$ as a renormalization factor in the motivic measure. In the equation (5) we deal with $\mathbb{L}$ as a class of the exceptional divisor minus two points. If we say that the "area" of this divisor may vary, we should add it as an additional parameter in the equations. 
Another motivation comes from the physical literature (e.g. \cite{phys}). There these parameters are interpreted as "Kaehler classes" of the exceptional divisors in the modification of $(\mathbb{C}^2,0)$.
It turns out that these parameters can be considered as independent times in a certain integrable system (\cite{phys}).

In any case, the process of solution of (\ref{lambda}) is completely analogous to the section 2.
We get $$G_{k,k}(t)=(\lambda_k-1)\sum_{m>k} G_{k,m}(t)=(\lambda_k-1)^2\sum_{m>k}t^{k(k-1)}\mathbb{L}^{-k}G_{k,m-k}(t)=$$

$$=(\lambda_k-1)^2\sum_{m>k}t^{k(k-1)}\mathbb{L}^{-k}(\sum_{m<k} G_{k,m}(t)+G_{k,k}(t)+\sum_{m>k} G_{k,m}(t))=$$

$$=(\lambda_k-1)^2\sum_{m>k}t^{k(k-1)}\mathbb{L}^{-k}(\sum_{m<k} G_{k,m}(t)+G_{k,k}(t)(1+\frac{1}{\lambda_k-1})),$$

so we get an analogue of the equation (6): 

\begin{equation}
G_{k,k}(t)=\frac{(\lambda_k-1)t^{k(k-1)}\mathbb{L}^{-k}}{1-\lambda_k\cdot t^{k(k-1)}\mathbb{L}^{-k}}\sum_{m<k} G_{k,m}(t).
\end{equation}

Since all equations except this one remains the same, we get the following analogue of the Theorem 2
(if we assume, that $G_{1,1}=(\mathbb{L}-1)^2\mathbb{L}^{-2}.$)

\begin{theorem}
Let $a=$g.c.d.$(k,m)$. Then
$$G_{k,m}(t,\mathbb{L})=$$
$$(\mathbb{L}-1)^2\cdot t^{(k-1)(m-1)}\mathbb{L}^{-k-m}\cdot \sum_{\underline{a},\underline{b}}\prod_{j=1}^{r}\mu(\frac{b_j}{a_{j-1}})\cdot \frac{(\lambda_{a_j}-1)t^{(a_j-1)b_j}\mathbb{L}^{-b_j}(1-t^{(a_j-1)(a_j-b_j)}\mathbb{L}^{b_j-a_j})}{(1-\lambda_{a_j}t^{a_j(a_j-1)}\mathbb{L}^{-a_j})\cdot (1-t^{(a_j-1)b_j}\mathbb{L}^{-b_j})},$$
where summation is indexed by all tuples $$(1=a_0\le b_1<a_1\ldots a_{r-1}\le b_r<a_r=a)$$ such that  $a_{j-1}|b_j$ and $b_j|a_j$
for all $j$ ( $r$ is not fixed).

\end{theorem}

\begin{corollary}
$$G_{k,m}(t^{-1},\mathbb{L}^{-1},\lambda_{1}^{-1},\lambda_{2}^{-1},\ldots)=t^{-2(k-1)(m-1)}\mathbb{L}^{2k+2m-2}\cdot G_{k,m}(t,\mathbb{L},\lambda_1,\lambda_2,\ldots).$$
\end{corollary}

\begin{proof}
This equation follows from the invariance of every of  $r$ multipliers in the right hand side of the theorem 1 under the change $(t,\mathbb{L},\lambda_{i})\mapsto (t^{-1},\mathbb{L}^{-1},\lambda_{i}^{-1}).$  
\end{proof}

Even for $t=1$ we get a nontrivial result.

\begin{corollary}
\begin{equation}
\label{t1}
G_{k,m}(1,\mathbb{L})=(\mathbb{L}-1)^2\mathbb{L}^{-k-m}\cdot \sum_{\underline{a},\underline{b}}\prod_{j=1}^{r}\mu(\frac{b_j}{a_{j-1}})\cdot \frac{(\lambda_{a_j}-1)\mathbb{L}^{-b_j}(1-\mathbb{L}^{b_j-a_j})}{(1-\lambda_{a_j}\mathbb{L}^{-a_j})\cdot (1-\mathbb{L}^{-b_j})},
\end{equation}
where summation is done as above.
\end{corollary}

Using parameters $\lambda_{k}$ we get a set of new identities that are hidden when all of them are equal to $\mathbb{L}$.

\begin{theorem}
Let $$H_{k,m}(\mathbb{L};\lambda_2,\lambda_3,\ldots)={G_{k,m}(1,\mathbb{L};\lambda_2,\lambda_3,\ldots)\over G_{k,m}(1,\mathbb{L};\mathbb{L},\mathbb{L},\ldots)}.$$
Define a generating function $$H(a,b,\mathbb{L};\lambda_2,\lambda_3,\ldots)=\sum_{k,m=1}^{\infty}a^{k}b^{m}H_{k,m}(\mathbb{L};\lambda_2,\lambda_3,\ldots).$$
Then the following equation holds:

\begin{equation}
\label{deq1}
{\partial \over \partial\lambda_{\alpha}}H(a,b,\mathbb{L};\lambda_2,\lambda_3,\ldots)=$$
$${1-\mathbb{L}^{-\alpha}\over (\lambda_{\alpha}-1)(1-\lambda_{\alpha}\mathbb{L}^{-\alpha})}\cdot H_{\alpha,\alpha}(\mathbb{L};\lambda_2,\lambda_3,\ldots,\lambda_{\alpha})\cdot H(a^{\alpha},b^{\alpha},\mathbb{L}^{\alpha};\lambda_{2\alpha},\lambda_{3\alpha},\ldots).
\end{equation}

More generally, if
$\lambda_{\alpha_1}<\lambda_{\alpha_2}<\ldots<\lambda_{\alpha_n},$ and $k_1,k_2,\ldots,k_n$ are some integers, then the higher derivatives can be expressed in a following form:

\begin{equation}
\label{deq2}
{\partial^{k_1+\ldots+k_n}\over \partial \lambda_{\alpha_1}^{k_1}\ldots \partial \lambda_{\alpha_n}^{k_n}}
H(a,b,\mathbb{L};\lambda_2,\lambda_3,\ldots)=$$
$$\prod_{j=1}^{n}[{(1-\mathbb{L}^{-\alpha_{j}})\mathbb{L}^{\alpha_{j}(k_j-1)}\over (\lambda_{\alpha_{j}}-1)(1-\lambda_{\alpha_j}\mathbb{L}^{-\alpha_{j}})^{k_j}}]\cdot H_{\alpha_1,\alpha_1}(\mathbb{L};\lambda_2,\lambda_3,\ldots)\times$$
$$\prod_{j=2}^{k}H_{{\alpha_{j}\over\alpha_{j-1}},{\alpha_{j}\over\alpha_{j-1}}}(\mathbb{L}^{\alpha_{j-1}};\lambda_{2\alpha_{j-1}},\lambda_{3\alpha_{j-1}},\ldots)\cdot H(a^{\alpha_n},b^{\alpha_n},\mathbb{L}^{\alpha_n};\lambda_{2\alpha_n},\lambda_{3\alpha_n},\ldots)
\end{equation}

if $\alpha_{j+1}\vdots \alpha_{j}$ for all $j=1\ldots n-1$, and otherwise the corresponding derivative vanishes.
\end{theorem}

\begin{proof}
Let $a=g.c.d.(k,m)$ and $\alpha|a$. We have
$$H_{k,m}=\sum_{\underline{a},\underline{b}}\prod_{j=1}^{r}\mu(\frac{b_j}{a_{j-1}})\cdot \frac{(\lambda_{a_j}-1)\mathbb{L}^{-b_j}(1-\mathbb{L}^{b_j-a_j})}{(1-\lambda_{a_j}\mathbb{L}^{-a_j})\cdot (1-\mathbb{L}^{-b_j})},$$
where $a_{j-1}\le b_j<a_j, a_{j-1}|b_j, b_j|a_j$. A summand in this sum has a nontrivial derivative by $\lambda_{\alpha}$, if $a_k=\alpha$ for some $k$. In this case we have a subsequence $1=a_0\le b_1<a_1\le b_2<\ldots<a_k=\alpha$ in every pair of sequences ${a},{b}$ which can be chosen independently of all subsequent $a_j$ and $b_j$, which are divisible by $\alpha$. Denote $\widetilde{a}_{j}={a_{k+j}\over \alpha},\widetilde{b}_{j}={b_{k+j}\over \alpha}$.
Therefore a the sum $S_{\alpha}$ of terms containing $\lambda_{\alpha}$ can be decomposed into a product
$$S_{\alpha}=H_{\alpha,\alpha}\cdot \sum_{\underline{\widetilde{a}},\underline{\widetilde{b}}}\prod_{j=1}^{r}\mu(\frac{\tilde{b}_j}{\tilde{a}_{j-1}})\cdot \frac{(\lambda_{\alpha\widetilde{a}_j}-1)\mathbb{L}^{-\alpha\widetilde{b}_j}(1-\mathbb{L}^{\alpha(\widetilde{b}_j-\widetilde{a}_j)})}{(1-\lambda_{\alpha\widetilde{a}_j}\mathbb{L}^{-\alpha\widetilde{a}_j})\cdot (1-\mathbb{L}^{-\alpha\widetilde{b}_j})},$$
where summation is indexed by all tuples $$(1=\widetilde{a}_0\le \widetilde{b}_1<\widetilde{a}_1\ldots \widetilde{a}_{r-1}\le \widetilde{b}_r<\widetilde{a}_r={a\over \alpha})$$ such that  $a_{j-1}|b_j$ and $b_j|a_j$.
This product can be rewritten as a product
$$S_{\alpha}=H_{\alpha,\alpha}(\mathbb{L};\lambda_2,\lambda_3,\ldots,\lambda_{\alpha})\cdot H_{{k\over\alpha},{m\over\alpha}}(\mathbb{L}^{\alpha};\lambda_{2\alpha},\lambda_{3\alpha},\ldots).$$

Now we remark that $${\partial H_{\alpha,\alpha}\over \partial\lambda_{\alpha}}={1-\mathbb{L}^{-\alpha}\over (\lambda_{\alpha}-1)(1-\lambda_{\alpha}\mathbb{L}^{-\alpha})}H_{\alpha,\alpha},$$
so $${\partial H_{k,m}\over\partial \lambda_{\alpha}}={1-\mathbb{L}^{-\alpha}\over (\lambda_{\alpha}-1)(1-\lambda_{\alpha}\mathbb{L}^{-\alpha})}H_{\alpha,\alpha}\cdot H_{{k\over\alpha},{m\over\alpha}}(\mathbb{L}^{\alpha};\lambda_{2\alpha},\lambda_{3\alpha},\ldots),$$
what implies the equation (\ref{deq1}).

The proof of the equation (\ref{deq2}) is completely analogous.
\end{proof}

We can try to specialize the result of the theorem 3 to some special choices of parameters $\lambda_i$.
For example, we can set $\lambda_{j}=A\tau^{j}$, where $A$ is a parameter and $\tau$ is a variable.

We get a function $Z(a,b,\mathbb{L},\tau)=H(a,b,\mathbb{L};A\tau^2,A\tau^3,\ldots)$.
From the equation (\ref{deq1}) we get immediately
\begin{equation}
{\partial\over \partial\tau}Z(a,b,\mathbb{L},\tau)=A\sum_{\alpha=2}^{\infty} \alpha\tau^{\alpha-1}{1-\mathbb{L}^{-\alpha}\over (A\tau^{\alpha}-1)(1-A\tau^{\alpha}\mathbb{L}^{-\alpha})}\cdot Z_{\alpha,\alpha}(\mathbb{L},\tau)\cdot Z(a^{\alpha},b^{\alpha},\mathbb{L}^{\alpha},\tau^{\alpha}).
\end{equation}

If we denote $Z_{k}(\mathbb{L},\tau)=H_{k,k}(\mathbb{L};A\tau^2,A\tau^3,\ldots),$
then we can also write a curious identity
\begin{equation}
{\partial\over\partial \tau}Z_n(\mathbb{L},\tau)=A\sum_{\alpha|n,\alpha\ge 2}\alpha\tau^{\alpha-1}{1-\mathbb{L}^{-\alpha}\over (A\tau^{\alpha}-1)(1-A\tau^{\alpha}\mathbb{L}^{-\alpha})}\cdot Z_{\alpha}(\mathbb{L},\tau)\cdot Z_{n\over\alpha}(\mathbb{L}^{\alpha},\tau^{\alpha}).
\end{equation}

Moscow State University,\newline
Department of Mathematics and Mechanics.

E.mail: gorsky@mccme.ru

\end{document}